\documentclass[reqno]{amsart} 

\usepackage{amsmath,amssymb}
\usepackage{amsthm}
\usepackage{indentfirst}
\usepackage{graphicx}
\usepackage{xcolor}
\usepackage[margin=1.5in]{geometry}
\usepackage{graphicx,tikz}

\theoremstyle{definition}

\theoremstyle{remark}

\begin{document}

\title[]{A nonlocal transport equation describing\\ roots of polynomials under differentiation}

\author[]{Stefan Steinerberger}
\address{Department of Mathematics, Yale University, New Haven, CT 06511, USA}
\email{stefan.steinerberger@yale.edu}

\thanks{ S.S. is partially supported by the NSF (DMS-1763179) and the Alfred P. Sloan Foundation.}

\begin{abstract} Let $p_n$ be a polynomial of degree $n$ having all its roots on the real line distributed according to a smooth function $u(0,x)$. One could wonder how the distribution
of roots behaves under iterated differentation of the function, i.e. how the density of roots of $p_n^{(k)}$ evolves. We derive a nonlinear
transport equation with nonlocal flux
$$ u_t + \frac{1}{\pi}\left( \arctan{ \left( \frac{Hu}{ u}\right)} \right)_x = 0 \qquad \mbox{on} ~\mbox{supp} \left\{u>0\right\},$$
where $H$ is the Hilbert transform. This equation has three very different compactly supported solutions: (1) the arcsine distribution $u(t,x) = (1-x^2)^{-1/2} \chi_{(-1,1)}$,
(2) the family of semicircle distributions
$$ u(t,x) = \frac{2}{\pi} \sqrt{(T-t) - x^2}$$
and (3) a family of solutions contained in the Marchenko-Pastur law.

\end{abstract}

\maketitle

\vspace{-15pt}

\section{Introduction}

\textbf{Introduction.} If $p_n$ is a polynomial of degree $n$ having $n$ distinct roots on the real line, then Rolle's theorem implies that $p_n^{(k)}$ has all its $n-k$ roots on the real line as well. 
Moreover, there is an interlacing phenomenon. A result commonly attributed to Riesz \cite{farmer, riesz} implies that the minimum gap between consecutive roots of $p_n'$ is bigger than that
of $p_n$: zeroes even out and become more regular. It is classical (and follows from interlacing) that if $p_n$ has its roots distributed according to
some nice distribution function, then $p_n'$ has its roots distributed according to the same function as $n \rightarrow \infty$.
The detailed study of the distribution of roots of $p_n'$ depending on $p_n$ is an active field \cite{bruj, bruj2, branko, dimitrov, gauss,  han,  kab, lucas, malamud, or, pem, rav, steini2, totik, ull, van}.
By the same reasoning, $p_n^{(k)}$ is also distributed following the same distribution for every fixed $k$ as $n \rightarrow \infty$.
However, this is no longer true when $k$ grows with $n$.
\begin{quote}
\textbf{Problem.} Let $(p_n)$ be a polynomial with $\mbox{deg} ~p_n = n$ and having only real roots whose distribution approximates a smooth probability distribution on $\mathbb{R}$ in a strong quantitative sense (say Kolmogorov-Smirnov or Wasserstein distance). What can be said about the
distribution of roots of $p_n^{(0.001n)}$? 
\end{quote}
If the roots of $p_n$ are evenly spaced at scale $\sim n^{-1}$, then interlacing implies that roots of derivative are shifted by at most $\sim n^{-1}$ which
implies that the dynamical evolution starts happening when the number of derivatives is comparable to the number of roots.\\

\textbf{The equation.} In the process of investigating this question, we came across a mean-field approximation that leads to a linear transport equation with nonlocal flux that can describe the evolution of the distribution of roots under iterated differentiation. The main purpose of this paper is to derive (in \S 3) and introduce the nonlinear equation
$$
\boxed{ u_t + \frac{1}{\pi} \left(\arctan{ \left( \frac{Hu}{ u}\right)} \right)_x = 0}
$$
where the equation is valid on the support $\mbox{supp}~u = \left\{x: u(x) > 0 \right\}$ and 
$$ Hf(x) =  \mbox{p.v.}\frac{1}{\pi} \int_{\mathbb{R}}{\frac{f(y)}{x-y} dy} \qquad \mbox{is the Hilbert transform.}$$
 The equation has the obvious symmetries under translation $u(x) \rightarrow u(x-\lambda)$ and reflection $u(x) \rightarrow u(-x)$. Moreover, since the Hilbert transform $Hu$ commutes with dilation, there
is an additional symmetry $u(t,x) \rightarrow \lambda u(t,x/\lambda)$. If $\mbox{supp}~u$ is an interval, then, assuming sufficient regularity ($u$ vanishing on the boundary of the support), we formally have
\begin{align*}
 \frac{d}{dt} \int_{\mathbb{R}}{u(x) dx} &= \int_{\mathbb{R}}{ u_t(x) dx} \\
 &= -\frac{1}{\pi} \int_{\mbox{\small supp}~u}{ \frac{d}{dx} \left(\arctan{ \left( \frac{Hu}{ u}\right)} \right) dx} = -1.
\end{align*}
This is in line with how the equation was derived: there should be a constant loss of mass since $p_n^{(t \cdot n)}$ has $(1-t)n$ roots. In particular, the solution should vanish in finite time at $t=1$.\\

\textbf{Related equations.} The equation is quite nonlinear but somewhat similar to a series of recently derived one-dimensional transport
equations with nonlocal flux given by the Hilbert transform or the fractional Laplacian. These were introduced as models for the quasi-geostrophic equation and one-dimensional analogues of the three-dimensional Navier-Stokes and Euler equations: we refer to Balodis \& Cordoba \cite{balodis}, Carrillo, Ferreira \& Precioso \cite{car},
Castro \& Cordoba \cite{castro},
Chae, Cordoba, Cordoba \& Fontelos \cite{cord1}, Constantin, Lax \& Majda \cite{const},
Cordoba, Cordoba \& Fontelos \cite{cord2}, Do, Hoang, Radosz \& Xu \cite{do},
 Dong \cite{dong}, Dong \& Li \cite{dong2},
Lazar \& Lemari\'{e}-Rieusset \cite{lazar}, Li \& Rodrigo \cite{li}
 and Silvestre \& Vicol \cite{silvestre}. Note added in print: Granero-Belinchon \cite{granero} has since studied an analogue of our equation on the one-dimensional torus.
We believe that it is conceivable that (a) techniques from that field could conceivably be useful in studying our transport equation (which is rather nonlinear) and (b) that,
conversely, the transport equation may be of interest in other contexts as well.

\section{Three explicit solutions}
We derive and describe three explicit compactly supported solutions in detail:
\begin{enumerate}
\item the stationary arcsine solution (not on all of $\mathbb{R}$ but only on $(-1,1)$)
$$u(t,x) = \frac{c}{\sqrt{1-x^2}} \chi_{(-1,1)} \qquad \mbox{where} \quad c \in \mathbb{R}$$
\item the semi-circle solution
$$ u(t,x) = \frac{2}{\pi} \sqrt{(T-t) - x^2} \cdot \chi_{|x| \leq \sqrt{T-t}} \qquad \mbox{for}~0 \leq t \leq T$$
\item the Marchenko-Pastur solution: introducing, for $c \geq 0$,
$$v(c,x)=  \frac{ \sqrt{(x_+ - x)(x - x_{-})}}{2 \pi x} \chi_{(x_{-}, x_+)}  \quad \mbox{where} \quad x_{\pm} = (\sqrt{c + 1} \pm 1)^2$$
that solution is given by
$$ u_c(t,x) =  v\left(\frac{c+t}{1-t}, \frac{x}{1 - t}\right).$$
\end{enumerate}

\subsection{The arcsine solution.}
We first describe the stationary solution when considering the equation only in the interval $(-1,1)$; in contrast to the other two solutions, the solution has singularities at the boundary of its support. If a function $f:\mathbb{R} \rightarrow \mathbb{R}$ is compactly supported on $(-1,1)$ and has its Hilbert transform $Hf$ vanish on its support, then it is
given by the arcsine distribution
$$u(t,x) = \frac{c}{\sqrt{1-x^2}} \chi_{(-1,1)} \qquad \mbox{where} \quad c \in \mathbb{R}.$$
This is true in a rather strong sense:
Coifman and the author \cite{coif} recently established that if $f(x)(1-x^2)^{1/4} \in L^2(-1,1)$ and $f(x) \sqrt{1-x^2}$ has mean value 0 on $(-1,1)$ (this enforces a form of orthogonality to the arcsine distribution), then
$$  \int_{-1}^{1}{ (Hf)(x)^2 \sqrt{1-x^2} dx} =  \int_{-1}^{1}{ f(x)^2 \sqrt{1-x^2} dx}.$$

This is mirrored in the classical fact that orthogonal polynomials on $(-1,1)$ with respect to a fairly large class of weights have their distribution of roots converge to the arcsine distribution
(see Erd\H{o}s \& Turan \cite{erd}, Erd\H{o}s \& Freud \cite{erd2}, Ullman \cite{ull} and Van Assche \cite{van}).

\begin{figure}[h!]
\begin{minipage}[l]{.45\textwidth}
\begin{tikzpicture}
\node at (0,0) {\includegraphics[width = 0.8\textwidth]{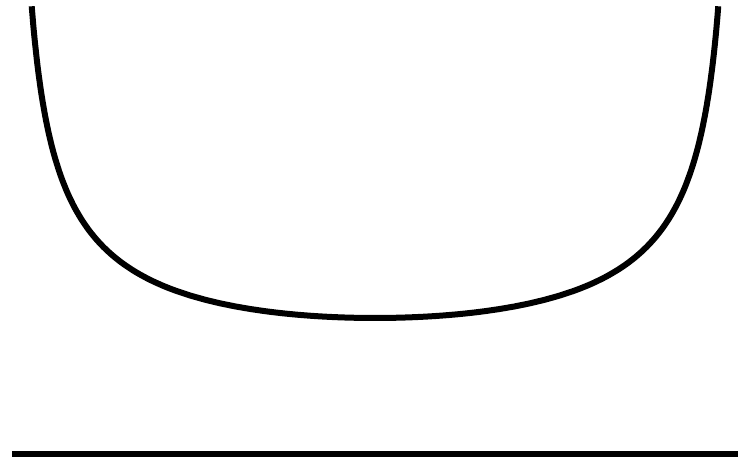}};
\end{tikzpicture}
\end{minipage} 
\begin{minipage}[r]{.45\textwidth}
\begin{tikzpicture}
\node at (0,0) {\includegraphics[width = 0.8\textwidth]{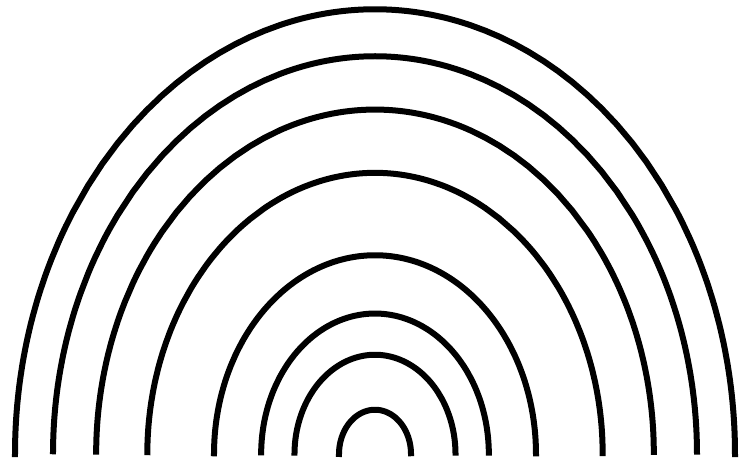}};
\end{tikzpicture}
\end{minipage} 
\caption{The arcsine solution (left) and the semicircle solution for $T=1$ (right) both shown for $t \in \left\{0,0.2, 0.4, 0.6, 0.8, 0.9, 0.95, 0.99\right\}$.} 
\end{figure}

\subsection{The semicircle distribution}
The construction of the semicircle solution is motivated by the behavior of the Hermite polynomials $H_n$. It is  known that
\begin{enumerate}
\item the roots of the Hermite polynomial $H_n$ are approximately (in the sense of weak convergence after rescaling) given by the measure
$$ \mu = \frac{1}{\pi} \sqrt{2n - x^2} dx$$
\item the derivatives of Hermite polynomials are again Hermite polynomials 
$$ \frac{d^m}{d x^m} H_n(x) = \frac{ 2^n n!}{(n-m)!} H_{n -m}(x).$$
\end{enumerate}
This suggests that if our transport equation models the flow of roots, then the semicircle solution should turn into a self-similar
one parameter family of solutions. A computation (carried out in \S 4.3) shows that, for every $T > 0$, 
$$ u(t,x) = \frac{2}{\pi} \sqrt{(T-t) - x^2} \cdot \chi_{|x| \leq \sqrt{T-t}} \qquad \mbox{for}~ t \leq T$$
is indeed a solution for $0 \leq t \leq T$.

\subsection{The Marchenko-Pastur solutions.} Our construction of the Marchenko-Pastur solution is motivated by the behavior of Laguerre polynomials. Laguerre polynomials $L_n$
do not form an Appell sequence, i.e. they are not closed under differentiation, however, the larger family of associated Laguerre polynomials $L_n^{(\alpha)}$ satisfies
$$ \frac{d^k}{dx^k} L_n^{(\alpha)}(x) = (-1)^k L_{n-k}^{(\alpha + k)}(x).$$
Moreover, the asymptotic distribution of roots is given by a Marchenko-Pastur distribution (indexed by a parameter $\alpha$): more precisely, it is classical \cite{kor} that for $n$ large, the roots
of $L_n^{(c \cdot n)}$ rescaled by a factor of $n$ converge in distribution to the Marchenko-Pastur distribution 
$$v(c,x)=  \frac{ \sqrt{(x_+ - x)(x - x_{-})}}{2 \pi x} \chi_{(x_{-}, x_+)} dx \quad \mbox{where} \quad x_{\pm} = (\sqrt{c + 1} \pm 1)^2.$$
Combining these two facts, we see that, asymptotically and for $ 0 < t < 1$,
$$ \frac{d^{t \cdot n}}{d x^{t \cdot n}} L_n^{c \cdot n} \sim \mbox{const} \cdot L_{(1-t) n}^{(c+t) \cdot n}$$
and this suggest that our nonlocal transport equation should have a solution of the form
$$ u_c(t,x) =  v\left(\frac{c+t}{1-t}, \frac{x}{1 - t}\right).$$
This is indeed the case.
For large values of $c$, the profile approximates that of the semicircle distribution (see Fig. 2). Presumably this will have implications for the stability analysis
around a semicircle distribution with Marchenko-Pastur solutions.
\begin{figure}[h!]
\begin{minipage}[l]{.45\textwidth}
\begin{tikzpicture}
\node at (0,0) {\includegraphics[width = \textwidth]{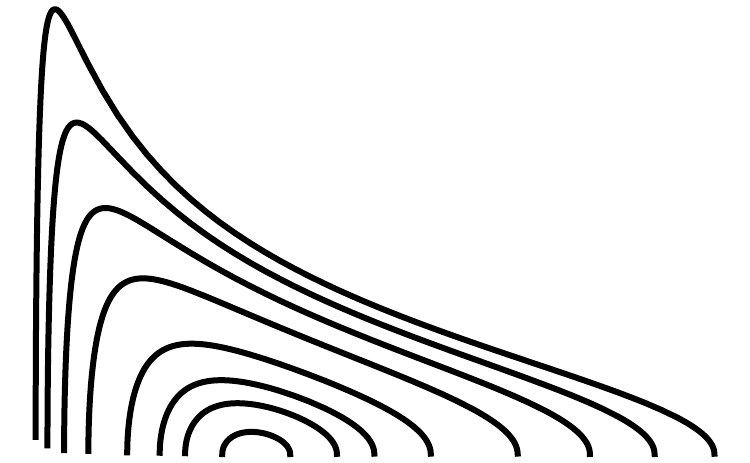}};
\end{tikzpicture}
\end{minipage} 
\begin{minipage}[r]{.5\textwidth}
\begin{tikzpicture}
\node at (0,0) {\includegraphics[width = \textwidth]{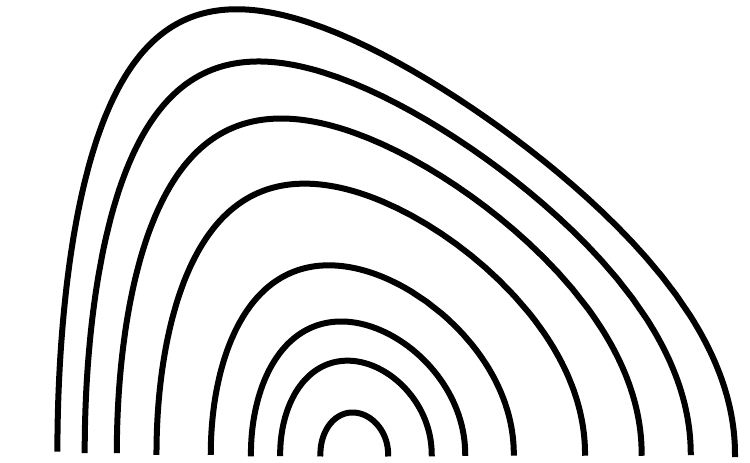}};
\end{tikzpicture}
\end{minipage} 
\caption{Marchenko-Pastur solutions $u_{c}(t,x)$: $c=1$ (left) and $c=15$ (right) shown for $t \in \left\{0,0.2, 0.4, 0.6, 0.8, 0.9, 0.95, 0.99\right\}$.} 
\end{figure}

\subsection{Outlook.} We believe that this motivates a rather large number of problems; it is natural to ask about the properties of the transport equation
itself: for which initial conditions is it well-posed? Is there a possibility of shock formation or finite-time blow-up? These questions might conceivably have
direct analogues for roots of polynomials under differentiation; presumably Riesz' theorem \cite{farmer, riesz} implies some basic form of regularity.
Another natural question is whether there is a rigorous derivation of the equation from polynomial dynamics in the small scale limit (this is likely to require a proper understanding of
the microstructure of roots). Are there other explicit solutions of the equation that that can be derived? What can be said about the stability properties
of the semicircle solution and the Marchenko-Pastur solution and does it correspond to polynomial dynamics? Finally, it seems natural to ask
whether there is an analogous equation (or possibly systems of equations) for polynomials with roots in the complex plane.

\section{Derivation of the equation}
\textbf{The Derivation.} Our derivation is based on two ingredients: (1) the Gauss interpretation of roots of derivatives as electrostatic equilibria (see \cite{gauss, lucas, steini})
and (2) Euler's cotangent identity. Regarding (1), we note that for any polynomial $p_n$ having roots in $\left\{x_1, \dots, x_n \right\} \subset \mathbb{R}$
$$ \frac{p_n'(x)}{p_n(x)} = \sum_{i=1}^{n}{\frac{1}{x-x_i}}.$$
This identity is also valid for polynomials in the complex plane with complex roots (thus suggesting that perhaps part of the derivation can be carried out
in the complex plane; what is missing is an analogue of the cotangent identity and the additional difficulty that density no longer uniquely defines a lattice).
The electrostatic interpretation also allows for an immediate proof of the Gauss-Lucas theorem \cite{gauss, lucas, steini, steini2}: the
roots of $p_n'$ are contained in the convex hull of the roots of $p_n$. This, in terms of our transport equation, implies that compactly supported initial
conditions give rise to compactly supported solutions (and that there is an inclusion relation for the support which is shrinking over time). Our second ingredient is the equation
$$ \pi \cot{\pi x} = \frac{1}{x} + \sum_{n=1}^{\infty}{ \left( \frac{1}{x+n} + \frac{1}{x-n} \right)} \qquad \mbox{for}~x \in \mathbb{R} \setminus \mathbb{Z}$$
dating back to Euler's \textit{Introductio in Analysis Infinitorum} (there is a particularly simple proof due to Herglotz \cite{bochner, herglotz}). We now assume that the roots of
a polynomial $p_n$ of very large degree $n$ are distributed according to a smooth density $u_0(x)$ and try to understand the microscopic movement of roots when passing
from $p_n$ to $p_n'$ at the local scale $n^{-1}$. Let us fix a root $p_n(y) = 0$. Recalling
$$ \frac{p_n'(x)}{p_n(x)} = \sum_{i=1}^{n}{\frac{1}{x-x_i}},$$
we split the right-hand side of that equation around $y$ into a far-field and a near-field. The far-field is approximately given by
$$  \sum_{|x_i - y| ~{\tiny \mbox{large}}}^{n}{\frac{1}{x-x_i}} \sim n\int_{\mathbb{R}}{\frac{u_0(y)}{x-y} dy} = n \pi (H u_0)(y),$$
where $H$ is the Hilbert transform. Here, $|x_i - y|$ being 'large' is to be understood as $n^{-1} \ll |x_i - y| \ll 1$. It remains to understand the near-field.
Since the distribution $u_0$ is smooth, the local density does not vary on short scales and we may approximate the near-field created by the local roots with a lattice structure; since the local
density is given by $u_0$, the spacing of the roots is given by $u_0(y)^{-1} n^{-1}$ and
\begin{align*}
  \sum_{|x_i - y| ~{\tiny \mbox{small}}}^{n}{\frac{1}{x-x_i}} &\sim \frac{1}{x-y_0} + \sum_{k \in \mathbb{N}}{ \left( \frac{1}{x - k u_0(y)^{-1} n^{-1}} + \frac{1}{x + k u_0(y)^{-1} n^{-1}} \right)}  \\
&= u_0(y) n \pi \cot{ (n \pi u_0(y) (x-y))}.
\end{align*}
The approximation is justified by the extremely fast convergence of the cotangent identity (assuming, of course, the underlying density to indeed be smooth).
Roots of $p_n'$ are created in places where the near-field and the far-field add up to 0, this leads to the equation
$$ u_0(y) \cot{ (n \pi u_0(y) (x-y))} =   (H u_0)(y).$$
This equation can be solved leading to
$$ x-y = -\frac{1}{n \pi} \arctan{\left( \frac{ (H u_0)(y)}{u_0(y)}\right)}$$
which informs us about the microscopix flux at scale $\sim n^{-1}$. This microscopic flow then gives rise to the transport equation
$$ u_t +\frac{1}{\pi}  \left(\arctan{ \left( \frac{Hu}{ u}\right)} \right)_x = 0.$$
\vspace{5pt}

\textbf{Missing ingredients.} There are two missing ingredients to making the derivation rigorous: (1) a rigorous understanding of the dynamics at the boundary of the support and (2) a rigorous understanding of what is happening in the bulk (this distinction is admittedly somewhat tautological).\\

(1) It is clear that our derivation, which assumes a flat background density of roots, must fail at the boundary where roots may have a different asymptotic behavior. We emphasize that ignoring these issues (which only affect a very small proportion of roots) seems to still result in a reasonable equation that is solved by at least three classical distributions. It might be that the contribution that the boundary has on the global dynamics is somewhat negligible but this remains to be rigorously proven.\\

(2) The derivation in the bulk is accurate as long as $u(\cdot, t)$ is essentially constant on length scales slightly larger than $n^{-1}$. This requires the equation to somehow undergo a smoothing effect: the spacing between the roots becomes more regular and the change in scale of spacing evens out. This would perhaps not be all that surprising: we refer to the paper of Farmer \& Rhoades \cite{farmer} discussing the possible existence of such a phenomenon and connecting it to a series of classical results (i.e. an argument of M. Riesz \cite{riesz} showing that the smallest gap between roots increases when going from $p_n$ to $p_n'$). These phenomena do not seem to be currently understood.

\section{Verification of the solutions}

\subsection{The arcsine solution.}
We recall an argument given by Coifman and the author in \cite{coif}: for this we introduce the Chebyshev polynomials $T_k$ (that will also play a role in the subsequent sections) via
$$T_{0}(x) =1, T_{1}(x) = x \quad \mbox{and} \quad T_{k+1}(x) = 2x T_{k}(x) - T_{k-1}(x).$$
as well as Chebyshev polynomials of the second kind $U_k$ given by 
$$U_{0}(x) =1, U_{1}(x) = 2x \quad \mbox{and} \quad U_{k+1}(x) = 2x U_{k}(x) - U_{k-1}(x).$$
These sequences of polynomials are orthogonal and for $n,m \geq 1$
$$ \frac{2}{\pi}\int_{-1}^{1}{ T_n(x) T_m(x) \frac{dx}{\sqrt{1-x^2}}} =  \delta_{nm}  \qquad \mbox{and} \qquad  \frac{2}{\pi}\int_{-1}^{1}{ U_n(x) U_m(x) \sqrt{1-x^2}dx} = \delta_{nm}.$$
The crucial identity is
$$   \frac{1}{\pi}\int_{-1}^{1}{\frac{a_k T_k(y)}{(x-y) \sqrt{1-y^2}} dy} = a_k U_{k-1}(x).$$
In particular, considering the function $g(x) = f(x) \sqrt{1-x^2}$ and expanding it into Chebyshev polynomials,
we see that the Hilbert transform acts as a shift operator. That shift operator annihilates exactly constants. The shift operator is also responsible for the fact that
 if $f(x)(1-x^2)^{1/4} \in L^2(-1,1)$ and $f(x) \sqrt{1-x^2}$ has mean value 0 on $(-1,1)$, then
$$  \int_{-1}^{1}{ (Hf)(x)^2 \sqrt{1-x^2} dx} =  \int_{-1}^{1}{ f(x)^2 \sqrt{1-x^2} dx}.$$
This shows that if $Hu$ vanishes on $(-1,1)$
for some $u$ compactly supported on $(-1,1)$, then $u$ is necessarily the arcsine distribution. This also shows that this is the only time-independent
solution of our transport equation when restricted to an open interval.

\subsection{The semicircle solution.}
As discussed above, the asymptotics of roots of Hermite polynomials combined with the fact that Hermite polynomials form an Appell sequence (closure under differentiation) suggests that
$$ u(t,x) = \frac{2}{\pi} \sqrt{(T-t - x^2)} \cdot \chi_{|x| \leq \sqrt{T-t}}\qquad \mbox{for}~0 \leq t \leq T$$
should be a solution of the equation. Clearly,
$$ \frac{\partial}{\partial t} u = -\frac{1}{\pi \sqrt{T-t-x^2}}.$$
It remains to compute the Hilbert transform $Hu$. The Hilbert transform commutes with positive dilations and is linear, we thus scale the function by a factor of $\sqrt{T-t}$ to reduce it to the computation
of the Hilbert transform of $(1-x^2)_{+}^{1/2}$ supported on $(-1,1)$. This reduces to a simple identity for Chebyshev polynomials of the second kind $U_k$
$$ \frac{1}{\pi}\int_{-1}^{1}{ \frac{\sqrt{1-y^2}U_{n-1}(y)}{x-y} dy} =  T_n (x)$$
since $U_0(x) = 1$ and $T_1(x)=x$ and thus, for $x$ in the support of $u$,
$$ H u(t,x) =\frac{2x}{\pi}  \chi_{(-\sqrt{T-t}, \sqrt{T-t})},$$
where $\chi$ is the characteristic function. A simple computation shows that
$$ \frac{1}{\pi}\left(\arctan{\left(\frac{Hu}{u}\right)}\right)_x = \frac{1}{\pi}\left(\arctan{ \frac{x}{\sqrt{(T-t - x^2)}}}\right)_x = \frac{1}{\pi \sqrt{T-t-x^2}}$$
and this shows that the semicircle solution solves the transport equation.

\subsection{The Marchenko-Pastur solution.}
Laguerre polynomials $L_n^{(\alpha)}$ are given by the recursion formula
$$ L_n^{(\alpha)}(x) =  \frac{x^{-\alpha}}{n!} \left( \frac{d}{dx} - 1 \right)^n x^{n+\alpha}.$$
Their behavior under differentiation is fairly easy to describe
$$ \frac{d^k}{dx^k} L_n^{(\alpha)}(x) = (-1)^k L_{n-k}^{(\alpha + k)}(x).$$
The behavior of the roots of $L_n^{(\alpha)}$ for $\alpha \geq 0$ is essentially classical \cite{bos,dette, gaw,kor,martinez}. The result that will inspire the construction of our solution uses
that if $\alpha_n$ is a sequence such that $\alpha_n/n \rightarrow c \in (-1, \infty)$, then the empirical distribution of the roots of $L_n^{(\alpha_n)}$
rescaled by a factor of $n$ converges weakly to the Marchenko-Pastur distribution
$$v(c,x)=  \frac{ \sqrt{(x_+ - x)(x - x_{-})}}{2 \pi x} \chi_{(x_{-}, x_+)} dx \quad \mbox{where} \quad x_{\pm} = (\sqrt{c + 1} \pm 1)^2.$$
Heuristically, we see that if
$$ \mbox{roots of}~L_{n}^{(c \cdot n)} \sim v(c,x) \quad \mbox{then} \quad \mbox{roots of}~L_{n(1-\varepsilon)}^{((c+\varepsilon) \cdot n)} \sim v\left(\frac{c+\varepsilon}{1-\varepsilon}, \frac{x}{1 - \varepsilon}\right)$$
This suggests the existence of a solution of the form
$$ u(t,x) =  v\left(\frac{c+t}{1-t}, \frac{x}{1 - t}\right).$$
We now verify the existence of the solution. The Hilbert transform of the Marchenko-Pastur distribution is known (see e.g. \cite[\S 5.5.2]{blower}) and given by
$$ H v(c,x) = \frac{x -c}{2\pi x} \qquad \mbox{on}~(x_{-}, x_{+}).$$
A  somewhat lengthy computation then shows that
$$ \frac{1}{\pi}\left(\arctan{\left(\frac{Hv(\frac{c+t}{1-t},\frac{1}{t-t})}{v(\frac{c+t}{1-t},\frac{1}{t-t})}\right)}\right)_x =  \frac{c+t+x}{2\pi  x \sqrt{2(2+c-t)x-(c+t)^2 - x^2}}$$
while
$$ \frac{\partial}{\partial t} v\left(\frac{c+t}{1-t}, \frac{x}{1 - t}\right) = - \frac{c+t+x}{2\pi  x \sqrt{2(2+c-t)x-(c+t)^2 - x^2}}$$
as desired.

\end{document}